\documentclass[10pt,a4paper]{article}
\usepackage{amsfonts}
\usepackage{amssymb}
\usepackage{mathrsfs}
\usepackage{amsthm}
\usepackage{setspace}
\usepackage{color}

\usepackage{fancyhdr}
\newtheorem{theo}{Theorem}[section]
\newtheorem{prop}[theo]{Proposition}

\theoremstyle{definition}

\newtheorem{defi}[theo]{Definition}

\newcommand{\be}{\begin{eqnarray*}}
\newcommand{\ee}{\end{eqnarray*}}
\newcommand{\beqa}{\begin{eqnarray}}
\newcommand{\eeqa}{\end{eqnarray}}
\newcommand{\ba}{\begin{array}}
\newcommand{\ea}{\end{array}}
\newcommand{\onab}{\overrightarrow{\nabla}}
\newcommand{\mc}{\mathcal}
\newcommand{\mf}{\mathfrak}
\newcommand{\ora}{\overrightarrow}
\newcommand{\rP}{\mathsf{P}}

\newcommand{\mbb}{\mathbb}

\begin{document}

\title{Note on pre-Courant algebroid structures for parabolic geometries}
\author{Stuart Armstrong}
\date{September 2007}
\maketitle

\begin{abstract}
This note aims to demonstrate that every parabolic geometry has a naturally defined per-Courant algebro\"id structure. If the geometry is regular, this structure is a Courant algebro\"id if and only if the the curvature $\kappa$ of the Cartan connection vanishes.
\end{abstract}

This note assumes familiarity with both parabolic geometry, and Courant algebro\"ids. See \cite{CartEquiv} and \cite{TCPG} for a good introduction to the first case, and \cite{Cdef} and \cite{TCA} for the second. Some of the basic definitions will be recalled here:
\begin{defi}[Parabolic Geometry]
Let $G$ be a semi-simple Lie group with Lie algebra $\mf{g}$, and $P$ a parabolic subgroup with Lie algebra $\mf{p}$. A parabolic geometry on a manifold $M$ is given by a principal $P$ bundle $\mc{P}$, an inclusion $\mc{P} \subset \mc{G}$, and a principal connection $\ora{\omega}$ on $\mc{G}$. This connection is required to satisfy the condition that $\ora{\omega}|_{\mc{P}}$ is a linear isomorphism $T\mc{P} \to \mf{g}$.
\end{defi}
Let $\mc{A} = \mc{P} \times_P \mf{g} = \mc{G} \times_G \mf{g}$ and denote by $\onab$ the affine connection on $\mc{A}$ coming from $\ora{\omega}$. By construction, $\mc{A}$ also inherits an algebra\"ic bracket $\{,\}$ and the Killing form $B$. It moreover has a well defined subbundle $\mc{A}_{(0)} = \mc{P} \times_P \mf{p}$, and the properties of $\onab$ give an equivalence $\mc{A} / \mc{A}_{(0)} \cong T$, thus a projection $\pi: \mc{A} \to T$. The Killing form $B$ then defines an inclusion $T^* \subset \mc{A}$, with $(T^*)^{\perp} = \mc{A}_{(0)}$. This implies that for $v$ a one form, $x$ any section of $\mc{A}$:
\be
B(v,x) = v \llcorner \pi(x).
\ee

Let $x,y$ be sections of $\mc{A}$. By the above inclusion, we may see $\onab x$ as section of $T^* \otimes \mc{A} \subset \mc{A} \otimes \mc{A}$, and thus we can directly write expressions like $\onab_x y$ (which is equal to $\onab_{\pi(x))} y$ in more traditional notation. The curvature of $\onab$ is $\kappa$; by inclusion, we can see this as a section of $\wedge^2 \mc{A} \otimes \mc{A}$. The parabolic structure gives a filtration on $\mc{A}$ and hence a concept of minimum homogeneity for sections of any associated bundle. There is also a Lie algebra co-differential $\partial^*: \wedge^2 \mf{p}^{\perp} \otimes \mf{g} \to \mf{p}^{\perp} \otimes \mf{g}$.
\begin{defi}
A parabolic geometry is \emph{regular} if $hom(\kappa) \geq 1$, and is normal if $\partial^* \kappa = 0$.
\end{defi}

Drawing on the definition of \cite{Cdef}:
\begin{defi}[Courant algebro\"id]
A Courant algebro\"id is vector bundle $E \to M$, with a pseudo-Riemannian metric $B$ on it, a projection $\pi: E \to TM$ called an anchor, and an inclusion $T^* \subset E$. It has a differential, $\mbb{R}$ -linear bracket $[,]: \Gamma(E) \otimes \Gamma(E) \to \Gamma(E)$. This is required to obey the following properties, for sections $x,y,z$ of $E$:
\begin{enumerate}
\item \label{pro1} The Jacobi identity $\mc{J}(x,y,z) = [x,[y,z]] - [[x,y],z] - [y,[x,z]] = 0$.
\item \label{pro2} $\pi(x) \cdot B(y,y) = 2 B(x,[y,y])$.
\item $\pi(x) \cdot B(y,z) = B([x,y],z) + B([x,z],y)$.
\end{enumerate}
\end{defi}
Note that property \ref{pro2} implies that $[,]$ is not a skew bracket. A pre-Courant algebro\"id, as defined in \cite{TCA}, is a structure that obeys all the pervious conditions, except for property \ref{pro1}. Instead, it is required to simply to have a linear Jacobian $\mc{J} \in \Gamma(\wedge^3 \mc{A} \otimes \mc{A})$.

The point of this note is:
\begin{theo}
Let $(M, \mc{A}, \onab)$ be a parabolic geometry. Then $\mc{A}$ is a pre-Courant algebro\"id for a natural choice of bracket $[,]$. If $\onab$ is flat, then it is a Courant algebro\"id. If the geometry is regular, then the Jacobian of $[,]$ has the same homogeneity as the curvature -- in particular, it is not a Courant algebro\"id whenever $\kappa \neq 0$.
into a Courant algebro\"id.
\end{theo}
First, it is easy to see that $\mc{A}$ has all the required algebra\"ic properties of a Courant algebro\"id: a projection to $T$, an inclusion of $T^*$, and a metric coming from the Killing form $B$.

Paper \cite{TCPG} defines a differential bracket $<,>$ defined on $\mc{A}$ as:
\be
<x,y> = \onab_x y - \onab_y x - \{x,y\} - \kappa(x,y).
\ee
Because of its original definition (defined as the bracket of right-invariant vector fields on $\mc{P}$), it must obey the Jacobi identity. This allows us to construct the (non-skew) bracket:
\be
[x,y] = <x,y> + B(y,\kappa(x,-)) - B(x,\kappa(y,-)) + B(\onab x, y).
\ee
The last term is the contraction of the second component of $\onab x$ with $y$; it is thus always a section of $T^*$. Most of the remaining properties of the pre-Courant algebro\"id are not hard to show (paper \cite{TCA} demonstrates them directly). For instance:
\beqa \label{eqa2}
[x,x] = B(\onab x, x) = \frac{1}{2} d B(x,x),
\eeqa
is immediate, (implying that $[x,y] = - [y,x] + d(B(x,y))$. On the other hand
\be
B([x,y],z) + B(y,[x,z]) &=& B(\onab_x y - \onab_y x - \{x,y\} - \kappa(x,y) + B(\onab x, y), z) \\
&& + B(y,\onab_x z - \onab_z x - \{x,z\} - \kappa(x,z) + B(\onab x, z)) \\
&& + B(B(y,\kappa(x,-),z) - B(B(x,\kappa(y,-),z) \\
&& + B(B(z,\kappa(x,-),y) - B(B(x,\kappa(z,-),y) \\
&=&- B(\{x,y\},z) - B(y,\{x,z\}) \\
&& + B(\onab_x y,z) - B(\onab_y x,z) + B(\onab_z x, y) \\
&& + B(\onab_x z,y) - B(\onab_z x,y) + B(\onab_y x, z) \\
&& - B(\kappa(x,y),z) - B(\kappa(x,z),y) \\
&& + B(y,\kappa(x,z)) - B(x,\kappa(y,z))\\
&& + B(z,\kappa(x,y)) - B(x,\kappa(z,y)) \\
&=& 0 + B(\onab_x y,z) + B(\onab_x z,y) \\
&=& \onab_x B(y,z) \\
&=& \pi(x) \cdot B(y,z),
\ee
since $\onab$ and $\{,\}$ preserve the metric $B$. This gives
\beqa \label{eqa1}
\pi(x) \cdot B(y,z) = B([x,y],z) + B(y,[x,z]).
\eeqa
Finally:
\begin{prop}
If $\mc{J} = [x,[y,z]] - [[x,y],z] - [y,[x,z]]$ is the Jacobian, then it is a section of $\wedge^3 \mc{A} \otimes \mc{A}$.
\end{prop}
\begin{proof}
First, we need to note that for any function $f \in C^{\infty}(M)$, $[df,y]=0$. This follows from the fact that $[df,y]$ is clearly a one-form (the only ambiguous term is $\onab_y df - \{y,df\}$, which is a one-form as the negative homogeneity components of $\onab_y$ and $\{y,-\}$ are the same) and from:
\be
[df,y](Z) &=& (- \onab_y df + \{y,df\} + B(\onab df, y) - B(df,\kappa(y,-)))(Z)\\
&=& - (\nabla df)(y,Z) - B(\{\rP(Y),df\}, Z) + B(\{A+v,df\},Z) + \\
&& (\nabla df)(Z,y) + B(\{Z,df\},y) +  B(\{\rP(Z),df\}, y) - B(df,\kappa(y,Z))\\
&=& d(df)(Z,y) - df(\nabla_Z Y - \nabla_Y Z -[Z,Y]) - B(df,\kappa(y,Z)) \\
&& + B(df, \{\rP(Y),Z\} -\{\rP(Z), y\} - \{A+v,Z\} + \{y,Z\}) \\
&=& B(df, (\{\rP(Y),Z\} -\{\rP(Z), Y\} + \{Y,Z\} + \nabla_Y Z - \nabla_Z Y -[Y,Z])_{-}) \\
&& - B(df,\kappa(y,Z))\\
&=& B(df, \kappa(Y,Z)_{-}) - B(df,\kappa(y,Z))\\
&=& 0,
\ee
where $Z$ is any section of $T$, and we have used a preferred connection $\nabla$ (see \cite{TCPG}) to give a splitting of $\onab_X = X + \nabla_X + \rP(X)$ and $y = Y + A + v$.

This allows us to calculate:
\be
\mc{J}(fx,y,z) &=& [fx,[y,z]] - [[fx,y],z] - [y,[fx,z]] \\
&=& f[x,[y,z]] -([y,z]\cdot f) x + B(x,[y,z])df \\
&& -[f[x,y],z] -[-(y \cdot f) x,z] -[B(x,y)df,z] \\
&& -[y,f[x,z]] +[y,(z \cdot f)x] - [y,B(z,x)df] \\
&=& f\mc{J}(x,y,z) -([y,z]\cdot f) x + B(x,[y,z])df \\
&& + (z\cdot f)[x,y] - B([x,y],z)df +(y \cdot f)[x,z] -(z\cdot y \cdot f)x + B(x,z)d(y\cdot f)\\
&&- B(x,y)[df,z] + (z \cdot B(x,y)) df - B(df,z) dB(x,y) \\
&&- (y\cdot f)[x,z] +(z\cdot f)[y,x] + (y\cdot z \cdot f) x - B(z,x)[y,df] - (y \cdot B(z,x)) df.
\ee
Re-arranging, and using the equations (\ref{eqa2}) and (\ref{eqa1}) extensively gives:
\be
&=& f\mc{J}(x,y,z) + B(x,z)(d(y\cdot f) -[y,df]) - B(x,y)[df,z] \\
&=& f\mc{J}(x,y,z) + B(x,z)[df,y] - B(x,y)[df,z] \\
&=& f\mc{J}(x,y,z).
\ee
Furthermore,
\be
\mc{J}(x,y,z) &=& [x,[y,z]] - [[x,y],z] - [y,[x,z]] \\
&=& -[x,[z,y]] + [z,[x,y]] + [[x,z],y] \\
&& +[x,dB(y,z)] - dB(z,[x,y]) - dB(y,[x,z]) \\
&=& - \mc{J}(x,z,y) + [x,dB(y,z)] - d(x\cdot B(y,z)) \\
&=& - \mc{J}(x,z,y) - [dB(y,z),x] \\
&=& - \mc{J}(x,z,y).
\ee
Similarly
\be
\mc{J}(x,y,z) &=& [x,[y,z]] - [[x,y],z] - [y,[x,z]] \\
&=& -[y,[x,z]] + ([[y,x],z] + [dB(y,x),z]) +[x,[y,z]]\\
&=& -\mc{J}(y,x,z).
\ee
So $\mc{J}$ is totally skew, and $C^{\infty}(M)$-linear in the first entry, hence in every entry.
\end{proof}

Given the previous result, to actually calculate $\mc{J}$ at a point, it suffices to choose a particular frame at that point. Let $\{e_j\}$ be a local frame of $\mc{A}$ around $p\in M$ chosen so that $(\onab e_j)_p = 0$. Then it is immediately evident that:
\be
[e_2,e_3]_p = (\{e_2,e_3\} + \kappa(e_2,e_3) + B(e_3,\kappa(e_2)) - B(e_2,\kappa(e_3)))_p.
\ee
Now consider $[e_1,[e_2,e_3]]_p$. The terms in that expression will be either second derivatives of the $e_j$, or linear terms. This gives us, at $p$:
\be
[e_1,[e_2,e_3]] &=& <e_1,<e_2,e_3>> +B(e_3, (\onab_{e_1} \kappa)(e_2)) - B(e_2, (\onab_{e_1} \kappa)(e_3)) \\
&& -\{e_1, B(e_3,\kappa(e_2)) - B(e_2,\kappa(e_3))\} + B(B(e_3,\kappa(e_2)) - B(e_2,\kappa(e_3)), \kappa(e_1)) \\
&& + B(\onab_{e_1} (\onab {e_2}), e_3)
\ee
\be
[e_2,[e_1,e_3]] &=& <e_2,<e_1,e_3>> +B(e_3, (\onab_{e_2} \kappa)(e_1)) - B(e_1, (\onab_{e_2} \kappa)(e_3)) \\
&& -\{e_2, B(e_3,\kappa(e_1)) - B(e_1,\kappa(e_3))\} + B(B(e_3,\kappa(e_1)) - B(e_1,\kappa(e_3)), \kappa(e_2)) \\
&& + B(\onab_{e_2} (\onab {e_1}), e_3)
\ee
\be
[[e_1,e_2],e_3] &=& - [e_3,[e_1,e_2]] + dB(e_3,[e_1,e_2]) \\
&=& -<e_3,<e_1,e_2>> -B(e_2, (\onab_{e_3} \kappa)(e_1)) + B(e_1, (\onab_{e_3} \kappa)(e_2)) \\
&& +\{e_3, B(e_2,\kappa(e_1)) - B(e_1,\kappa(e_2))\} - B(B(e_2,\kappa(e_1)) - B(e_1,\kappa(e_2)), \kappa(e_3)) \\
&& - B(\onab_{e_3} (\onab {e_1}), e_2) \\
&& + B(e_3, \onab (\onab_{e_1} e_2 - \onab_{e_2} e_1) - (\onab \kappa)(e_1,e_2) + B(e_2,(\onab \kappa)(e_1)) - B(e_1,(\onab \kappa)(e_2)) \\
&& + B(e_3,B(\onab(\onab e_1), e_2))
\ee
Now $<,>$ obeys the Jacobi identity. We get further simplifications of the type $(\onab_{e_1} \kappa)(e_2) = (\onab_{e_2} \kappa)(e_1) -(\onab \kappa)(e_1,e_2)$ by the Bianci identity on $\onab$. Remembering the identity $B(e_3,v) = v(e_3)$ for any one-form $v$ gives simplifications in the $[[e_1,e_2],e_3]$ term. Together, these give the relation:
\be
\mc{J}(e_1,e_2,e_3) &=& -\{e_1, B(e_3,\kappa(e_2)) - B(e_2,\kappa(e_3))\} + B(e_3,\kappa(e_2,\kappa(e_1))) - B(e_2,\kappa(e_3,\kappa(e_1)))  \\
&& + \textrm{ cyclic terms}.
\ee
This demonstrates that if $\kappa = 0$, then we are in the presence of a Courant algebro\"id. However, if $\kappa \neq 0$ and $\onab$ is regular, this construction will always fail to be a Courant algebro\"id:
\begin{prop}
If the Cartan connection is regular and $\kappa \neq 0$ at a point then $\mc{J} \neq 0$ at that point, and $hom(\mc{J}) = hom(\kappa)$.
\end{prop}
\begin{proof}
Assume $\kappa \neq 0$ at $p$. Let $h = hom(\kappa)$. Since $\kappa$ is regular, $h>0$. Since $B$ and $\{,\}$ preserve homogeneity,
\be
hom(\mc{J}) \geq h.
\ee
At $p$, let us project $\kappa$ onto its lowest homogeneity component $\kappa_H$ (if $\kappa$ is normal, this is just the lowest homogeneity harmonic curvature \cite{two}). This $\kappa_H$ can further be decomposed into sections $\kappa_{a,b,c}$ of $T^*_a \wedge T^*_b \otimes \mc{A}_c$ for integers $a,b$ and $c$, with $a + b +c = h$. Pick $a,b$ and $c$ such that $\kappa_{a,b,c} \neq 0$ at $p$.

The terms in $\mc{J}$ with two appearances of $\kappa$ are of homogeneity $\geq 2h > h$, so we will ignore them. Chose a Weyl structure $\nabla$ that preserve a volume form to give a splitting $\mc{A} = \sum_{i=-k}^k \mc{A}_i$. The Killing form $B$ ensures that $\mc{A}_i \perp \mc{A}_j$ whenever $j \neq -i$. Define $\mc{A}_{(j)} = \sum_{i=j}^k \mc{A}_i$ (this does not depend on the choice of $\nabla$).

Let $E$ be the grading section in $\mc{A}_0$, $X_{-a}$ a section of $\mc{A}_{-a} = T_{-a}$ and $Y_{-c}$ a section of $\mc{A}_{-c}$. Call $\mc{J}_h$ the homogeneity $h$ component of $\mc{J}$. Now $\kappa_H(E) = 0$, and assume for the moment that $c \neq 0$, implying that $B(E,\kappa_H(X_{-a}))$ has trivial projection onto the $T_{b}^*$ factor. We will now look, in homogeneity $h$, at the $T_{b}^*$ factor of $\mc{J}(E,X_{-a},Y_{-c})$:
\be
&\Big(\{E, B(Y_{-c}, \kappa_H(X_{-a}))\} - \{E, B(X_{-a}, \kappa_H(Y_{-c}))\} \Big) \\
&= \\
&\{E, (Id-m)(\kappa_H)(X_{-a}, - , Y_{-c}) \},
\ee
where $m$ is the operator interchanging the first and last component of $\otimes^3 \mc{A}$. Basic representation theory implies that $Id - m$ is injective on $\wedge^2 \mc{A} \otimes \mc{A}$, hence on $\kappa$. Moreover $Id - m$ preserves homogeneity, so there must exist $X_{-a}$, $Y_{-c}$ such that $(Id-m)(\kappa_H)(X_{-a}, - , Y_{-c})/(\mc{A}_{(b+1)})$ is a non-zero section of $\mc{A}_{b}$ around $p$. The bracket with $E$ does not change this, as $E$ acts by multiplication by $b$ on these sections, and $b>0$ since $\kappa$ is a curvature.

Now if $c = 0$, then let $Z_0$ be a section of $\mc{A}_0$, and then
\be
(\mc{J})_b(X_{-a}, Y_0, Z_0) = (Id-m)(\phi \circ \kappa_H(X_{-a})) (Y_0,Z_0).
\ee
Here $\phi$ is the map $\mc{A} \to \wedge^2 \mc{A}$ given by the Lie bracket. Since $\mc{A}$ is semi-simple, $\phi$ has no kernel, making $\phi \circ \kappa_H(X_{-a})$ into an non-degenerate section of $\wedge^2 \mc{A} \otimes \mc{A}$. But $Id-m$ is injective on this bundle, ensuring that there must exist $Y_0$ and $Z_0$ making the above expression non-zero.

This demonstrates that $\mc{J} \neq 0$ whenever $\kappa \neq 0$ and further, that $hom(\mc{J}) = h$.

\end{proof}

\bibliographystyle{alpha}
\bibliography{ref}

\begin{thebibliography}{{\v{C}}ap06}

\bibitem[{\v{C}}ap06]{two}
Andreas {\v{C}}ap.
\newblock Two constructions with parabolic geometries.
\newblock {\em Rend. Circ. Mat. Palermo (2)}, 79:11--37, 2006.

\bibitem[{\v{C}}G02]{TCPG}
Andreas {\v{C}}ap and Rod Gover.
\newblock Tractor calculi for parabolic geometries.
\newblock {\em Trans. Amer. Math. Soc.}, 354(4):1511--1548, 2002.

\bibitem[{\v{C}}S00]{CartEquiv}
Andreas {\v{C}}ap and Hermann Schichl.
\newblock Parabolic geometries and canonical {C}artan connections.
\newblock {\em Hokkaido Math. J.}, 29(3):453--505, 2000.

\bibitem[KS05]{Cdef}
Yvette Kosmann-Schwarzbach.
\newblock Quasi, twisted, and all that$\ldots$in {P}oisson geometry and {L}ie
  algebroid theory.
\newblock {\em Progr. Math.}, 232:363--389, 2005.

\bibitem[Vai05]{TCA}
Izu Vaisman.
\newblock Transitive {C}ourant algebroids.
\newblock {\em Int. J. Math. Math. Sci.}, 11:1737--1758, 2005.

\end{thebibliography}

\end{document}